\documentclass[11pt]{amsart}

\usepackage{mathrsfs}
\usepackage{amssymb,amsmath,amsxtra,amscd, feynmf,latexsym}
\usepackage{pstricks}
\usepackage{mathpple}
\usepackage{setspace,ifpdf}
\usepackage{diagrams}

\usepackage{pst-tree, pst-node}
\usepackage[mathscr]{eucal} 
\ifpdf
\usepackage{pdfsync}
\fi

\theoremstyle{plain}
\numberwithin{equation}{section}
\newtheorem{theorem}{Theorem}

\numberwithin{lemma}{section}
\numberwithin{corollary}{section}
\numberwithin{proposition}{section}

\theoremstyle{definition}
\newtheorem{definition}{Definition}
\numberwithin{definition}{section}

\theoremstyle{plain}
\newtheorem{example}{Example}
\numberwithin{example}{section}
\newtheorem{remark}{Remark}

\newcommand{\nc}{\newcommand}
\nc{\C}{\mathcal{C}}
\nc{\CC}{\widetilde{C}}
\nc{\wt}{\overline}
\nc{\mc}{\mathcal}
\nc{\on}{\operatorname}
\nc{\Cl}{Cl}
\nc{\drva}{\widehat{\OOmega}}
\nc{\spec}{\OOn{Spec}}
\nc{\AutO}{\OOn{Aut} \mc{O}}
\nc{\vac}{|0\rangle}
\nc{\Z}{\mathbb{Z}}
\nc{\zf}[1]{z^{\frac{1}{#1}}}
\nc{\wf}[1]{w^{\frac{1}{#1}}}
\nc{\Mt}{M^{\sigma}}

\nc{\gr}{\OOn{gr}}
\nc{\T}{\mathbb{T}}
\nc{\LT}{\mathbb{LT}}
\nc{\CT}{\mathbb{C}\{ \mathbb{T} \}}
\nc{\g}{\mathfrak{g}}
\nc{\A}{\mathcal{A}}
\nc{\n}{\mathfrak{n}}
\nc{\I}{\mathcal{I}}
\nc{\Hc}{\mathcal{H}}
\nc{\HH}{\mathbf{H}}
\nc{\M}{M}
\nc{\U}{\mathcal{U}}
\nc{\F}{\mathcal{F}}
\nc{\OO}{\mathcal{O}}
\nc{\LRF}{\mathcal{LRF}}
\nc{\FD}{\mathcal{LFG}}
\nc{\LFD}{\mathcal{LFG}}
\nc{\FG}{\mathcal{LFG}}
\nc{\LFG}{\mathcal{LFG}}
\nc{\labb}{\emph{lab}}

\nc{\Iso}{\on{Iso}}

\nc{\heck}{\mc{H}ecke}
\renewcommand{\F}{\mathcal{F}}

\nc{\al}{\alpha}
\nc{\OOl}{\OOverline}

\begin{document}

\title{Rooted trees, Feynman graphs, and Hecke correspondences}
\author{Matt Szczesny}
\address{Department of Mathematics  
         Boston University, Boston MA, USA}
\email{szczesny@math.bu.edu}

%\date{June 2009}

\begin{abstract}

We construct natural representations of the Connes-Kreimer Lie algebras on rooted trees/Feynman graphs arising from Hecke correspondences in the categories $\LRF, \LFG$ constructed by K. Kremnizer and the author. 
We thus obtain the insertion/elimination representations constructed by Connes-Kreimer as well as an isomorphic pair we term top-insertion/top-elimination. We also construct graded finite-dimensional sub/quotient representations of these arising from "truncated" correspondences. 

\end{abstract}

\maketitle 

\section{Introduction}

The Connes-Kreimer Hopf algebras on rooted trees and Feynman graphs  $\mc{H}_{\T}, \mc{H}_{FG}$, introduced in \cite{K}, \cite{CK3}, describe the algebraic structure of the BPHZ algorithm in the renormalization of perturbative quantum field theories. If we let $\T$ denote the set of (non-planar) rooted trees, and $\mathbb{Q}\{ \T \}$ the $\mathbb{Q}$--vector space spanned by these, then as an algebra, $\Hc_{\T} = \on{Sym}(\mathbb{Q}\{  \T \})$, and the coalgebra structure is given by the coproduct
\[
\Delta(T) = \sum_{C \textrm{ admissible cut }} P_C(T) \otimes R_{C} (T)
\]  
where $P_C(T)$ is the forest of branches resulting from the cut $C$, and $R_C(T)$ is the root component remaining "above" the cut (see \cite{CK3} for a more detailed definition). 

$\Hc_{FG}$ is defined analogously, with Feynman graphs in place of rooted trees. More precisely, given a perturbative QFT, and denoting by $\mathbb{Q}\{ \Gamma \}$ the vector space spanned by the one-piece irreducible graphs of the theory ($1\on{PI}$ graphs), $\Hc_{FG} = \on{Sym}(\mathbb{Q}\{ \Gamma \})$ as an algebra. Its coalgebra structure is given by 
\[
\Delta(\Gamma) = \sum_{\gamma \in \Gamma} \gamma \otimes \Gamma / \gamma
\]
where the sum is over all (not necessarily connected) subgraphs of $\Gamma$, and $\Gamma / \gamma$ denotes the graph obtained from $\Gamma$ by shrinking each connected component of $\gamma$ to a point. 

$\Hc_{T}$ and $\Hc_{FG}$ are graded connected commutative Hopf algebras, and so by the Milnor-Moore theorem, their duals $\Hc^{*}_{T}$ and $\Hc^{*}_{FG}$ are isomorphic to the universal enveloping algebras $\U(\n_T)$, $\U(\n_{FG})$ of the nilpotent Lie algebras $\n_T$, $\n_{FG}$ of their primitive elements. We refer to $\n_T$ and $\n_{FG}$ as the Connes-Kreimer Lie algebras on rooted trees and Feynman graphs respectively. 

In \cite{KS}, a categorification of the Hopf algebras $\U(\n_T)$, $\U(\n_{FG})$ was obtained, by showing that they arise naturally as the \emph{Ringel-Hall} algebras of certain categories $\LRF$, $\LFG$ of labeled rooted forests and Feynman graphs respectively. We briefly recall this notion. Given an abelian category $\C$, linear over a finite field, and having the property that all $\on{Hom}$'s and $\on{Ext}^1$ are finite-dimensional (such a category is called \emph{finitary}), one can construct from it a Hopf algebra $\on{H}_{\C}$, called its Ringel-Hall algebra. As a vector space, $\on{H}_{\C}$ is the space of $\mathbb{Q}$--valued functions on $\on{Iso}(\C)$ - the set of isomorphism classes of $\C$,  with finite support. The algebra structure is given by the associative product
\[
f \times g (M) := \sum_{N \subset M} f(N)g(M/N) \hspace{1cm} f, g \in \on{H}_{\C}, \; \; M \in \on{Iso} (\C),
\]
and the co-algebra structure by the co-commutative coproduct 
\begin{align*}
\Delta: \on{H}_{\C} & \rightarrow \on{H}_{\C} \otimes \on{H}_{\C} \\
\Delta(f)(M, N) & := f(M \oplus N)
\end{align*}
$\on{H}_{\C}$ is a graded connected co-commutative Hopf algebra, and so the enveloping algebra of the Lie algebra $\n_{\C}$ of its primitive elements. In this light, the main result of \cite{KS} is the construction of categories $\LRF, \LFG$ such that $\n_{\LRF} \simeq \n_{\T}$, and $\n_{\LFG} \simeq \n_{FG}$ (note however that $\LRF, \LFG$ are not abelian). 

Ringel-Hall algebras provide a very useful perspective in the study of the representations of $\n_{\C}$. In particular, we may construct representations of $\n_{\C}$ via \emph{Hecke correspondences}. Let 
\[
\heck_{\C} := \{ (A, B)  | B \in \on{Iso}(\C), A \subset B \}
\]
(i.e. the set of pairs consisting of an isomorphism class $B$ in $\C$, and a subobject $A$
of $B$). 
$\heck_{\C}$ comes with three maps to $\on{Iso}(\C)$ : $\pi_1, \pi_2, \textrm{ and } \pi_q$, where $\pi_1(A,B) = A$, $\pi_2 (A,B ) = B$, and $\pi_q (A, B) = B/A$. As $\on{H}_{\C}$ is the space of functions on $\on{Iso}(\C)$, we obtain the diagram
\begin{diagram}
   &   &   \F(\heck_{\C}) &  &    \\
   & \ruTo<{\pi^*_q} & \uTo^{\pi^*_1} & \luTo>{\pi^*_2} & \\
 \on{H}_{\C} &   &  \on{H}_{\C}   &   &   \on{H}_{\C} \\  
\end{diagram}
where $ \F(\heck_{\C})$ denotes the space of functions on $\heck_{\C}$. Since the fibers of  $\pi_1, \pi_2, \pi_q$ are finite, we also have maps $\pi_{1*}, \pi_{2*}, \pi_{q*}$ corresponding to integration along the fiber:
\[
\pi_{r*} (f) (y) = \sum_{x_y \in \pi_r^{-1}(y)} f(x_y) \hspace{1cm} r=1,2,q
\]
By taking $r,s,t$ to be some permutation of $1,2,q$, we now obtain convolutions:
\begin{align*}
\pi_{r,s,t} : \on{H}_{C} \otimes \on{H}_{\C} & \rightarrow \on{H}_{\C} \\
f \otimes g & \rightarrow \pi_{r*}(\pi^*_{s} ( f ) \pi^*_{t} ( g ))
\end{align*}
$\pi_{r,s,t}$ can be viewed as a map $$\on{H}_{\C} \rightarrow \on{Hom}(\on{H}_{\C}, \on{H}_{\C}).$$
and provided that $\pi_{r,s,t}$ is compatible with the algebra structure on $\on{H}_{\C}$, this yields a representation of $\n_{\C}$ on $\on{H}_{\C}$. Representations constructed in this manner are naturally graded by $K_{0}(\C)$. One checks easily that among them are the left and right actions of $\on{H}_{\C}$ on itself as well as their graded duals. 

The utility of this perspective extends far beyond recovering the action of $\n_{\C}$ on its enveloping algebra however. We may consider various sub-correpondences of $\heck_{\C}$. For example, we may choose an isomorphism class $M \in \on{Iso}(\C)$, and consider
 \[
\heck_{\leq M} := \{ (A,B) | A \subset B \subset M \}
\]
Convolution now yields an action of $\on{H}_{\C}$ on the finite-dimensional space $\on{H}_{\leq M}$ of functions supported on isomorphism classes of sub-objects of $M$. One may similarly obtain representations in $\on{H}_{quot(M}$--functions supported on isomorphism classes of quotient objects of $M$. 

These techniques are well-known in the case when $\C$ is a finitary abelian category, such as the category of modules over a finite-dimensional algebra. In this paper, we apply them in the case when $\C$  is one of the categories $\LRF, \LFG$, to study the representation theory of the Connes-Kreimer Lie algebras $\n_{T}, \n_{FG}$. These categories are not abelian, but share many properties of finitary abelian categories.  
We show that the insertion/elimination actions of $\n_{\C}$ introduced in \cite{CK} arise in this way, as do two "new" (but isomorphic to the previous two) representations of $\n_{\C}$, which we call "top insertion" and "top elimination". While the elimination and top-elimination (resp. insertion and top-insertion) representations are isomorphic, they have non-isomorphic finite-dimensional sub/quotient representations arising from correspondences of the form $\heck_{\leq M}, \heck_{quot(M)}$ as above, and have a combinatorially distinct flavor. One of the advantages of this approach is that it allows us to treat both the cases of rooted trees and Feynman graphs in one language. 

This paper is structured as follows. In sections \ref{Hall_alg}, \ref{LRF}, \ref{LFG}, and\ref{gen_properties} we recall the notion of Ringel-Hall algebra, the construction of categories $\LRF, \LFG$, and their general properties. Section \ref{left_right} gives explicit combinatorial formulas for the left/right actions of $\n_{\C}$ on $\on{H}_{\C} \simeq U(\n_{\C})$ and its graded dual. Sections \ref{fdrep1} and \ref{fdrep2} examine finite-dimensional graded representations derived from these. 
 In section \ref{Hecke_corr} we study Hecke correspondences in $\LRF, \LFG$, and show that we recover the representations from \ref{left_right}. We also show how to use "truncated" correspondences to construct the finite-dimensional representations from \ref{fdrep1} and \ref{fdrep2}.  \bigskip

\noindent{\bf Acknowledgements:} I would like to thank Dirk Kreimer, Valerio Toledano-Laredo, and Kobi Kremnizer for many valuable conversations.  I would also like to thank Olivier Schiffmann for suggesting the use of Hecke convolutions in the representation theory of $\n_{\T}, \n_{FG}$. 

\section{Ringel-Hall algebras}
\label{Hall_alg}
This section briefly recalls the notion of Ringel-Hall algebra, following \cite{S, BT}. 
Recall that a small abelian category $\C$ is called \emph{finitary} if: 
\begin{align}
\textrm{ i) For any two objects } M,N \in \on{Ob}(\C) \textrm{ we have } |Hom(M,N)| < \infty \\
\textrm{ii) For any two objects } M,N \in \on{Ob}(\C) \textrm{ we have } |Ext^1 (M,N)| < \infty
\end{align}

Denote by $\on{Iso}(\C)$ the set of isomorphism classes of objects in $\C$. 

\begin{definition}
The Ringel-Hall algebra of $\C$, denoted $\on{H}_{\C}$ is the vector space of finitely supported $\mathbb{Q}$--valued functions on $\on{Iso}(\C)$. I.e.
\[
\on{H}_{\C} := \{ f: \on{Iso}(\C) \rightarrow \mathbb{Q} | |supp(f)| < \infty \}
\]
\end{definition}

\noindent $\on{H}_{\C}$ is equipped with an associative product defined by
\[
f \times g (M) := \sum_{N \subset M} f(N)g(M/N) \hspace{1cm} f, g \in \on{H}_{\C}, \; M \in \on{Iso} (\C).
\]
and a co-commutative coproduct 
\begin{align*}
\Delta: \on{H}_{\C} & \rightarrow \on{H}_{\C} \otimes \on{H}_{\C} \\
\Delta(f)(M, N) & := f(M \oplus N)
\end{align*}

The two structures are easily seen to be compatible, and $\on{H}_{\C}$ is in fact a co-commutative Hopf algebra, graded by $K_{0} (\C)$ (via the natural map $\on{Iso}(\C) \rightarrow K_0 (\C)$), and connected.  By the Milnor-Moore theorem, it is the enveloping algebra ot the Lie algebra $\n_{\C}$ of its primitive elements, called the \emph{Ringel-Hall Lie algebra} of $\C$. 

In the applications considered in this paper, the category $\C$ in question will not be abelian, but "nearly" so, and the Ringel-Hall algebra construction goes through without problems (see \cite{KS} for a detailed discussion of these issues). 

\section{The category $\LRF$}
\label{LRF}
In this section, we briefly recall the construction of the category $\LRF$ of labeled rooted forests from \cite{KS}. We begin by reviewing some notions related to rooted trees. Let $S$ be an infinite set. For a tree $T$, denote by $V(T), E(T)$ the vertex and edge sets of $T$ respectively.   

\begin{definition}
\begin{enumerate}
\item A \emph{rooted tree labeled by $S$} is a tree $T$, with a distinguished vertex $r(T)\in V(T)$, and an injection $l: V(T) \hookrightarrow S$. Denote by $RT(S)$ the set of all such. 
\item A \emph{rooted forest labeled by $S$} is a set $F$ of rooted trees labeled by $S$, whose labels are disjoint, i.e. $$F=\{T_1, T_2, \cdots, T_k  \}, \; \; T_i \in RT(S), im(l_i) \cap im(l_j) = \emptyset \textrm{ if } i \neq j.$$
\item An \emph{admissible cut} of a labeled tree $T$ is a subset $C(T) \subset E(T)$ such that at most one member of $C(T)$ is encountered along any path joining a leaf to the root. Removing the edges in an admissible cut divides $T$ into a labeled rooted forest $P_C(T)$ and a labeled rooted tree $R_C(T)$, where the latter is the component containing the root. The \emph{empty} and \emph{full} cuts $C_{null}, C_{full}$, where $$(P_{C_{null}}(T), R_{C_{null}}(T)) = (\emptyset, T) \textrm{ and }  (P_{C_{full}}(T), R_{C_{full}}(T)) = (T, \emptyset)$$ respectively, are considered admissible.
\item An \emph{admissible cut} on a labeled forest $F = \{ T_1, \cdots, T_k \}$ is a collection of cuts  $C=\{C_1, \cdots, C_k \}$, with $C_i$ an admissible cut on $T_i$. Let
\begin{align*}
R_C (F) & := \{ R_{C_1}(T_1), \cdots, R_{C_k}(T_k) \} \\
P_C (F) & := P_{C_1}(T_1) \cup P_{C_2} (T_2) \cup \cdots \cup P_{C_k} (T_k)
\end{align*}
\item The \emph{maximum} of two admissible cuts $C,D$ on a rooted tree $T$, denoted by $max(C,D)$, is the admissible cut obtained by taking the cut edge closer to the root along any path from leaf to root. Similarly, we define $min(C,D)$ by taking the cut edge further from the root along any path from leaf to root. 
\item Two labeled rooted forests $F_1= \{ T_1, \cdots, T_k \}$ and $F_2 = \{ T'_1, \cdots, T'_m \}$ are isomorphic if $k=m$, and there is a permutation $\sigma \in S_k$ and bijections $$f_i : V(T_i) \rightarrow V(T'_{\sigma(i)}),  \; \; i=1, \cdots ,k. $$ which preserve roots and all other incidences. 
\end{enumerate}
\end{definition}

\begin{example} Consider the labeled rooted forest consisting of a single tree $T$, with root drawn at the top. 

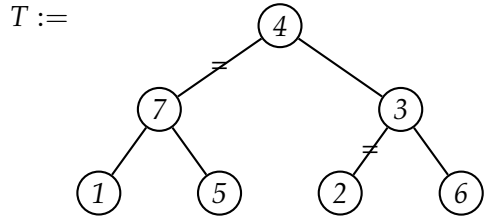
\begin{center} \psset{levelsep=6ex, treesep=1.0cm}
$T:=$ \pstree{  \Tcircle{4}  } {  \pstree{ \Tcircle{7} \ncput{=}   } { \Tcircle{1} \Tcircle{5}}
 \pstree{\Tcircle{3}}{\Tcircle{2}\ncput{=} \Tcircle{6}}   }
\end{center}
and the cut edges are indicated with "=", then 
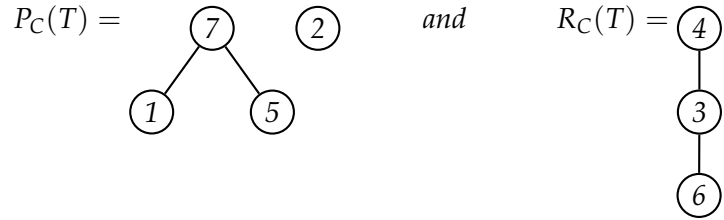
\begin{center}  \psset{levelsep=6ex, treesep=1.0cm}
$P_C(T) = $ \pstree{ \Tcircle{7}}{\Tcircle{1} \Tcircle{5}} \hspace{3cm} \Tcircle{2} \hspace{1cm}
and
\hspace{1cm}
$R_C(T) = $ \pstree{\Tcircle{4}} {\pstree{\Tcircle{3}} {\Tcircle{6}} }
\end{center}

\end{example}

\noindent We are now ready to define the category $\LRF$. 

\begin{definition} The category $\LRF$ is defined as follows:
\begin{itemize}
\item
 \[
\on{Ob}(\LRF) = \{ \textrm{ labeled rooted forests } \} \cup \{ \emptyset \}
\]
where $\emptyset$ denotes the \emph{empty forest}, which plays the role of zero object.
\item 
\begin{align*}
\on{Hom}(F_1,F_2) := & \{ (C_1,C_2,f) | C_i \textrm{ is an admissible cut of } F_i, \\ & \; f: R_{C_1} (F_1) \cong P_{C_2} (F_2) \} \; \; F_i \in \on{Ob}(\LRF). 
\end{align*}
( For $F \in \LFG$, $(C_{null}, C_{full}, id): F \rightarrow F$ is the identity morphism in $\on{Hom}(F,F)$. ) 
\end{itemize}

\end{definition}

\bigskip

\noindent {\bf Example:} if 
\bigskip
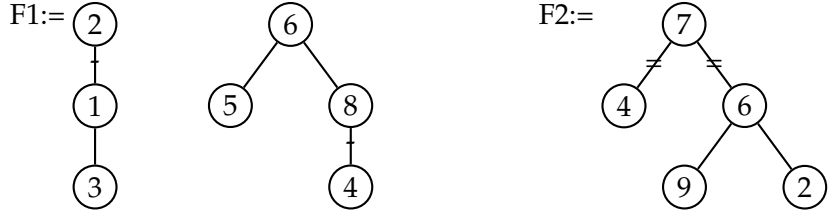
\begin{center} \psset{levelsep=6ex, treesep=1.0cm}
F1:= \pstree[]{ \Tcircle{2} }{ \pstree{ \Tcircle{1} \ncput{-} } {\Tcircle{3} } } \hspace{1 cm} 
\pstree{\Tcircle{6}}{\Tcircle{5} {\pstree{\Tcircle{8}}{\Tcircle{4}\ncput{-} } }} \hspace{2cm}
F2:= \pstree{\Tcircle{7}}{ \Tcircle{4}\ncput{=} \pstree{\Tcircle{6}\ncput{=} }{\Tcircle{9}\Tcircle{2}}}
\end{center}
\bigskip
then a morphism is given by the triple $(C_1, C_2,f)$ where:
\begin{itemize}
\item $C_1$ is indicated by "$-$", and $C_2$ is indicated by "$=$". 
\item $f : R_{C_1}(F_1) \cong P_{C_2}(F_2)$ is defined by $f(2) = 4$, $f(5) = 9$, $f(6)=6$ ,$f(8)=2$. 
\end{itemize}

\noindent We recall the definition of the the composition of morphisms
\[
\on{Hom}(F_1, F_2) \times \on{Hom}(F_2, F_3) \rightarrow \on{Hom}(F_1,F_3)
\]
Suppose that $(C_1,C_2,f) \in \on{Hom}(F_1,F_2)$, and $(D_2,D_3,g) \in \on{Hom}(F_2,F_3)$. The cut $min(C_2,D_2)$ induces a cut $E_1$ on $F_1$, and $max(C_2,D_2)$ a cut $E_3$ on $F_3$. The restriction of $g \circ f$ gives an isomorphism  $R_{E_1}(F_1) \cong P_{E_3} (F_3)$. We define the composition above to be $(E_1,E_3, g \circ f)$. The associativity of composition follows from the associativity of $max$ and $min$.

\subsection{The Connes-Kreimer Lie algebra on rooted trees}

In this section, we recall the definition of the Connes-Kreimer Lie algebra on rooted trees $\n_T$ (see \cite{CK}). As a vector space, 
\[
\n_T = \mathbb{Q}\{ \T \} 
\]
i.e. the span of unlabeled rooted trees. On $\n_T$, we have a \emph{pre-Lie} product
"$*$", given, for $T_1, T_2 \in \T$ by
\[
T_1 * T_2 = \sum_{T \in \T} a(T_1,T_2;T) T
\]
where 
\[
a(T_1,T_2;T) := | \{ e \in E(T) | P_{C_e}(T) = T_1, R_{C_e} (T) = T_2 \} |
\]
and $C_e$  denotes the cut severing the edge $e$. 
The Lie bracket on $\n_T$ is given by 
\begin{equation} \label{treebracket}
[T_1,T_2] := T_1 * T_2 - T_2 * T_1
\end{equation}
Thus, for example if
%\begin{align*}
%\psset{levelsep=0.3cm, treesep=0.3cm}
%[\pstree{\Tr{\bullet}}, \pstree{\Tr{\bullet}}{\Tr{\bullet}\Tr{\bullet}} ]&= \pstree{\Tr{\bullet}}{\Tr{\bullet} 
%\Tr{\bullet} \Tr{\bullet }} + \pstree{\Tr{\bullet}}{\pstree{\Tr{\bullet}}{\Tr{\bullet}}\Tr{\bullet}}         - \pstree{\Tr{\bullet}}{\pstree{\Tr{\bullet}}{ \Tr{\bullet} \Tr{\bullet}} } 
%\end{align*}
\begin{center} \psset{levelsep=4ex, treesep=1.0cm}
$T_1:=$ \pstree{\Tcircle{}} {\Tcircle{} {\Tcircle{}} }  \hspace{2cm} \textrm{ and } $T_2:=$ \pstree{\Tcircle{}}{}
\end{center}
then 
\begin{center} \psset{levelsep=4ex, treesep=1.0cm}
$[T_1, T_2] = $ \pstree{\Tcircle{}} {\pstree{\Tcircle{}} {\Tcircle{} \Tcircle{}}} $-$  \pstree{\Tcircle{}}{ {\pstree{\Tcircle{}}{\Tcircle{}}} \Tcircle{}} $- 3 $ \pstree{\Tcircle{}}{\Tcircle{} \Tcircle{} \Tcircle{}}
\end{center}

\section{The category $\LFG$}
\label{LFG}
In this section we review the construction of the category $\LFG$ of labeled Feynman graphs following \cite{KS}. Our treatment of the combinatorics of graphs is taken from \cite{Y}. In order to not get bogged down in notation, we focus on the special case of $\phi^3$ theory (the case of trivalent graphs with only one edge-type). The results of this section extend to the general case in a completely straighforward manner.

\begin{definition}
A \emph{graph} $\Gamma$ consists of a set $H = H(\Gamma)$ of half-edges, a set $V = V(\Gamma)$ of vertices, a set of vertex-half edge adjacency relations $( \subset V \times H)$, and a set of half edge - half edge adjacency relations $(\subset H \times H)$, with the requirements that each half edge is adjacent to at most one other half edge and to exactly one vertex. Note that graphs may not be connected. 

Half edges which are not adjacent to another half edge are called \emph{external edges}, and denoted $Ex=Ex(\Gamma) \subset E = E(\Gamma)$. Pairs of adjacent half edges are called \emph{internal edges}, and denoted $Int(\Gamma)$. 
\end{definition}

\begin{definition}
A \emph{half edge $S$--labeled graph}, (\emph{labeled graph} for short), is a triple $(\Gamma, S, \rho)$, where $\Gamma$ is a graph, $S$ is a set such that $|S|=|H|$, and  $\rho: H \rightarrow  S $ is a bijection. $S$ will usually be obvious from context. 
\end{definition}

\begin{definition}
\begin{enumerate}
\item A \emph{Feynman graph} is a graph where each vertex is incident to exactly three half-edges, each connected component has $2$ or $3$ external edges, and has at least one loop (i.e. $H_1 (\Gamma) \geq 1)$.  We denote the set of Feynman graphs by $FG$. 

\item Similarly, we can define  a (half-edge) \emph{labeled Feynman graph}. We denote the set of labeled Feynman graphs by $LFG$.

\item We say that two labeled Feynman graphs $\Gamma_1$ and $\Gamma_2$ are \emph{isomorphic} if there exist bijections $f_V: V(\Gamma_1) \rightarrow V(\Gamma_2)$, $f_H: H(\Gamma_1) \rightarrow H(\Gamma_2)$ which induce bijections on all incidences. We write $f: \Gamma_1 \cong \Gamma_2$.

\item A graph (or a labeled graph) is \emph{1-particle irreducible} ($1\on{PI}$) if it is connected, and remains connected under the removal of an arbitrary internal edge. 
\end{enumerate}
\end{definition}

\noindent {\bf Example:} The graph $\Gamma_{eg}$

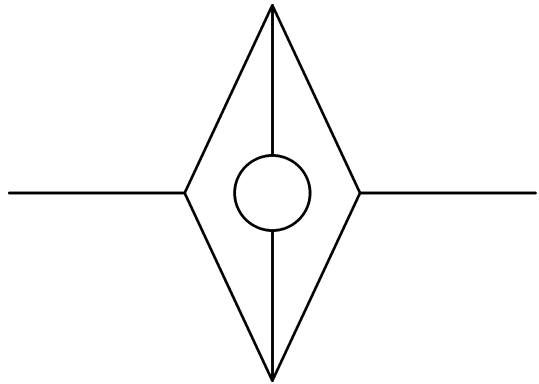
\begin{center} \label{graph1}
\unitlength=1mm
\begin{fmffile}{fig1} 
\begin{fmfgraph*}(70,50)
\fmfleft{i}
\fmfright{o}
\fmftop{t}
\fmfbottom{b}
\fmflabel{$v5$}{b}
\fmflabel{$v6$}{t}
\fmflabel{$v1$}{v1}
\fmflabel{$v2$}{v2}
\fmflabel{$v3$}{v3}
\fmflabel{$v4$}{v4}
\fmf{plain, label=$e1$}{i,v1}
\fmf{plain, label=$e2$}{v2,o}
\fmf{plain, label=$e3$}{t,v1}
\fmf{plain, label=$e4$}{t,v2}
\fmf{plain, label=$e5$}{b,v1}
\fmf{plain, label=$e6$}{b,v2}
\fmf{plain, label=$e7$}{b,v3}
\fmf{plain, label=$e8$}{t,v4}
% \fmf{plain,left,tension=1}{v3,v4,v3}
\fmf{plain, left, label=$e9$, tension=1}{v3,v4}
\fmf{plain, right, label=$e10$, tension=1}{v3,v4}
\end{fmfgraph*}
\end{fmffile}
\end{center}
\vspace{1cm}

Is a $1\on{PI}$ Feynman graph with two external edges. We have labeled each vertex and edge, and each half-edge can be thought of as labeled by a pair $(v,e)$ where $v$ is a vertex, and $e$ is an edge incident to $v$. 

\begin{definition}
Given a Feynman graph $\Gamma$, a \emph{subgraph} $\gamma$ is a Feynman graph such that  $V(\gamma) \subset V(\Gamma)$, $H(\gamma) \subset H(\Gamma)$, and such that if $v \in V(\gamma)$, and $(v,e) \in V(\Gamma) \times H(\Gamma)$, then $e \in H(\gamma)$ (i.e. the subgraph has to contain all half-edges incident to its vertices). We write $\gamma \subset \Gamma$.  The same definition applies to labeled graphs. 
\end{definition}
\bigskip
\noindent {\bf Example:} We define a subgraph $\gamma_{eg} \subset \Gamma_{eg}$ as follows. Let $V(\gamma_{eg}) = \{ v3, v4 \}$, $E(\gamma_{eg}) := \{ \textrm{ all half-edges incident to } v3, v4 \}$, i.e. $$E(\gamma_{eg}) = \{ (v4,e8), (v3,e7), (v3,e9), (v3,e10), (v4,e9), (v4,e10)  \}$$ and all incidences inherited from $\Gamma_{eg}$. 
\bigskip

\noindent We proceed to define the contraction of subgraphs of Feynman graphs. 

\begin{definition}
Let $\Gamma$ be a Feynman graph, and $\gamma \subset \Gamma$ a connected subgraph. The \emph{quotient graph} $\Gamma / \gamma$ is defined as follows. If $\gamma$ has $3$ external edges, then $\Gamma / \gamma$ is the Feynman graph with 
\begin{enumerate}
\item $V(\Gamma/ \gamma)$ set the vertex set of $\Gamma$ with all vertices of $\gamma$ removed, and a new trivalent vertex $v$ added. 
\item $H(\Gamma/ \gamma)$ the half edge set of $\Gamma$, with all half edges corresponding to internal half edges of $\gamma$ removed.
\item All adjacencies inherited from $\Gamma$, and the external half edges of $\gamma$ joined to $v$.  
\end{enumerate}
\bigskip
If $\gamma$ has $2$ external edges, then $\Gamma/ \gamma$ is the Feynman graph with
\bigskip
\begin{enumerate}
\item $V(\Gamma / \gamma)$ is $V(\Gamma)$ with all the vertices of $\gamma$ removed.
\item $H(\Gamma / \gamma)$ is $H(\Gamma)$ with all half edges of $\gamma$ removed. 
\item All adjacencies inherited from $\Gamma$, as well as the adjacency of the external half-edges of $\gamma$. 
\end{enumerate}
\bigskip
Finally, If $\gamma \subset \Gamma$ is an arbitrary (not necessarily connected) Feynman subgraph, then $\Gamma / \gamma$ is defined to be the Feynman graph obtained by performing successive quotients by each connected component. Note that the order of collapsing does not matter. 
\end{definition}
\bigskip

\noindent {\bf Example:} With $\Gamma_{eg}, \gamma_{eg} $ as above, $\Gamma_{eg} / \gamma_{eg}$ is:
\vspace{1cm}
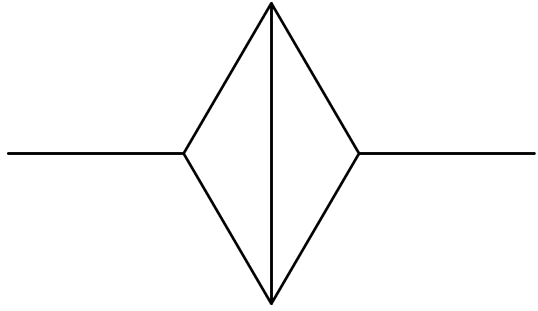
\begin{center} \label{fig2}
\unitlength=1mm
\begin{fmffile}{fig2} 
\begin{fmfgraph*}(70,40)
\fmfleft{i}
\fmfright{o}
\fmftop{t}
\fmfbottom{b}
\fmflabel{$v5$}{b}
\fmflabel{$v6$}{t}
\fmflabel{$v1$}{v1}
\fmflabel{$v2$}{v2}
\fmf{plain, label=$e1$}{i,v1}
\fmf{plain, label=$e2$}{v2,o}
\fmf{plain, label=$e3$}{t,v1}
\fmf{plain, label=$e4$}{t,v2}
\fmf{plain, label=$e5$}{b,v1}
\fmf{plain, label=$e6$}{b,v2}
\fmf{plain, label=$e7$}{t,b}
\end{fmfgraph*}
\end{fmffile}
\end{center}
\vspace{1cm}

\begin{remark} \label{remark3}
If $\gamma_1, \gamma_2 \subset \Gamma$ are subgraphs of a (labeled) Feynman graph, then $\gamma_1 \cap \gamma_2$ is a Feynman subgraph of $\gamma_i$ (and of course also $\Gamma$).
\end{remark}
\bigskip
\begin{remark} \label{remark4}
If $\gamma \subset \Gamma$ is a subgraph of a (labeled) Feynman graph, then there is a bijection between subgraphs of $\Gamma/ \gamma$, and subgraphs $\gamma'$ of $\Gamma$ such that $\gamma \subset \gamma' \subset \Gamma$.
\end{remark}
\bigskip

\noindent Labeled Feynman graphs form a category $\LFG$ as follows:

\begin{itemize}
\item \[
\on{Ob}(\LFG) = \{ \textrm{ labeled Feynman graphs } \} \cup \{ \emptyset \}
\]
where $\emptyset$ denotes the \emph{empty Feynman graph}, which plays the role of zero object. Note that objects of $\LFG$ may have several connected components. 
\item If $\Gamma_1, \Gamma_2 \in \LFG$, we define
\begin{align*}
\on{Hom}(\Gamma_1, \Gamma_2) := & \{ (\gamma_1, \gamma_2, f) | \gamma_i \textrm{ is a subgraph of } \Gamma_i, \\ & f: \Gamma_1 / \gamma_1 \cong \gamma_2 \}
\end{align*}
(For $\Gamma \in \LFG$, $(\emptyset,\Gamma,id)$ is the identity map in $\on{Hom}(\Gamma,\Gamma)$.)

\end{itemize}

\noindent The composition of morphisms in $\LFG$
\[
\on{Hom}(\Gamma_1, \Gamma_2) \times \on{Hom}(\Gamma_2, \Gamma_3) \rightarrow \on{Hom}(\Gamma_1,\Gamma_3)
\]
is defined as follows. Suppose that $(\gamma_1,\gamma_2,f) \in \on{Hom}(\Gamma_1,\Gamma_2)$, and $(\tau_2,\tau_3,g) \in \on{Hom}(\Gamma_2,\Gamma_3)$. By remark \ref{remark3}, the subgraph $\tau_2$ on $\Gamma_2$ induces a subgraph $\tau_2 \cap \gamma_2$ of $\gamma_2 \subset \Gamma_2$, which by remark \ref{remark4} corresponds to a subgraph of $\xi$ of $\Gamma_1$ containing $\gamma_1$. The image $g \circ f (\xi) \subset \Gamma_3$ is a subgraph $\rho \subset \tau_3$. We define the composition $(\tau_2,\tau_3,g) \circ (\gamma_1,\gamma_2,f)$ to be $(\xi,\rho, g\circ f)$. It is easy to see that this composition is associative.

\subsection{The Connes-Kreimer Lie algebra on Feynman graphs}

In order to define the Connes-Kreimer Lie algebra structure on Feynman graphs, we must first introduce the notion of inserting a graph into another graph. Let
$$\mathbb{Q} \{ LFG \} $$ denote the vector space spanned by labeled Feynman graphs. 

\begin{definition}
Let $\Gamma_1, \Gamma_2 \in LFG$. If $\Gamma_1$ has three external edges, $v \in V(\Gamma_2) $, and $f: Ex(\Gamma_1) \rightarrow H(v)$ is a bijection (where $H(v)$ are the labeled half-edges incident to the vertex $v$), then let $\Gamma_2 \circ_{v,f} \Gamma_1$ be the labeled Feynman graph such that
\bigskip
\begin{itemize}
\item $V(\Gamma_2 \circ_{v,f} \Gamma_1) = V(\Gamma_2) \cup V(\Gamma_1) \backslash v$.
\item $H(\Gamma_2 \circ_{v,f} \Gamma_1) = H(\Gamma_1) \cup_f H(\Gamma_2)$ - i.e. the unions of the half-edges of each graph, with the identifications induced by $f$.
\item The adjacencies induced from those of $\Gamma_1$ and $\Gamma_2$. 
\end{itemize}
\bigskip
If $\Gamma_1$ has two external edges, $\{ e_1,e_2 \} \in Int(\Gamma_2) \subset H \times H$, and $f$ is a bijection between $Ex(\Gamma_1)$ and $\{ e_1, e_2 \}$ (there are two of these), then $\Gamma_2 \circ_{e,f} \Gamma_1$ is the labeled Feynman graph such that
\bigskip 
\begin{itemize}
\item $V(\Gamma_2 \circ_{e,f} \Gamma_1) = V(\Gamma_1) \cup V(\Gamma_2)$.
\item $H(\Gamma_2 \circ_{e,f} \Gamma_1) = H(\Gamma_1) \cup H(\Gamma_2)$.
\item The adjacency induced by $f$ as well as those induced from $\Gamma_1$ and $\Gamma_2$.
\end{itemize}
\bigskip

Let $\n_{LFG}$ denote the $\mathbb{Q}$--vector space spanned by \emph{unlabeled} connected Feynman graphs. Given a labeled Feynman graph $\Gamma$, denote by $\wt{\Gamma}$ the corresponding unlabeled Feynman graph. Thus, 
\[
\n_{LFG} = \mathbb{Q} \{ LFG \} / \sim
\]
where $\Gamma \sim \Gamma' $ iff $\wt{\Gamma} = \wt{\Gamma'}$.
\noindent We now equip $\n_{LFG}$ with the pre-Lie product "$\star$", defined by
\[
\wt{\Gamma_1} \star \wt{\Gamma_2} := \sum_{v \in V(\Gamma_2), \\
 f: Ex(\Gamma_1) \rightarrow H(v)}  \wt{\Gamma_2 \circ_{v,f} \Gamma_1 }
\]
if $\Gamma_1$ has three external edges, and
\[
\wt{\Gamma_1} \star \wt{\Gamma_2} := \sum_{e \in Int(\Gamma_2), \\ f: Ex(\Gamma_1) \rightarrow \{ e_1, e_2 \} } \wt{\Gamma_2 \circ_{v,f} \Gamma_1}
\]
if $\Gamma_1$ has two external edges, and extended linearly ( in the above formulas, we first choose an arbitrary labeling of the Feynman graphs).  
 Finally, we can define the Lie bracket on $\n_{LFG}$ by
\begin{equation} \label{FGLiebracket}
[\wt{\Gamma_1},\wt{\Gamma_2}] := \wt{\Gamma_1} \star \wt{ \Gamma_2 } - \wt{\Gamma_2} \star \wt{ \Gamma_1}
\end{equation}

\end{definition}  
\bigskip
\begin{fmffile}{fig6}
\unitlength=1mm
\noindent {\bf Example:} Suppose
\begin{align}
\Gamma_1 &=
\parbox{20mm}{
\begin{fmfgraph}(20,20)
\fmfleft{i1}
\fmfright{o1}
\fmf{plain}{i1,v1}
\fmf{plain}{v2,o1}
\fmf{plain,left,tension=0.7}{v1,v2,v1}
\end{fmfgraph}} &
\Gamma_2 &= 
\parbox{20mm}{
\begin{fmfgraph}(40,40)
\fmfleft{i}
\fmftop{t}
\fmfbottom{b}
\fmf{plain}{i,v1}
\fmf{plain}{t,v2}
\fmf{plain}{b,v3}
\fmf{plain}{v1,v2,v3,v1}
\end{fmfgraph}}
\end{align}
\end{fmffile}

\noindent then 

\unitlength=1mm
\begin{fmffile}{fig7}
\begin{align*}
[\Gamma_1,\Gamma_2] &= 6
\parbox{20mm}{
\begin{fmfgraph}(40,40)
\fmfleft{i}
\fmftop{t}
\fmfbottom{b}
\fmf{plain}{i,v1}
\fmf{plain}{t,v2}
\fmf{plain}{b,v3}
\fmf{plain}{v2,v3}
\fmf{plain}{v3,v1}
\fmf{plain}{v1,v4}
\fmf{plain}{v2,v5}
\fmf{plain,left,tension=0.8}{v4,v5,v4}
\end{fmfgraph}}
& - 12 \hspace{1cm} &
\parbox{20mm}{\begin{fmfgraph}(40,40)
\fmfleft{i}
\fmfright{o}
\fmftop{t}
\fmfbottom{b}
\fmf{plain}{i,v1}
\fmf{plain}{v2,o}
\fmf{phantom}{t,v3}
\fmf{phantom}{b,v4}
\fmf{plain}{v1,v3,v2,v4,v1}
\fmf{plain}{v3,v4}
\end{fmfgraph}}
\end{align*}
\end{fmffile}

\section{Properties of the categories $\LRF, \LFG$} \label{properties}
\label{gen_properties}
In this section, we summarize some of the properties of $\LRF, \LFG$ that allow us to define the Ringel-Hall algebra of these categories. For more detail, please consult\cite{KS}. 

Let $\C$ denote either of the categories $\LRF, \LFG$. Then $\C$ has the following properties:

\begin{enumerate}
\item The empty forest/graph $\{ \emptyset \}$ is an null object. 
\item Disjoint union of forests/graphs equips $\C$ with a symmetric monoidal structure, which we denote by $\oplus$. We will refer to objects of the form $A \oplus B$, with neither $A$ nor $B$ the null object $\emptyset$ as \emph{decomposable}, and those not admitting such a splitting as \emph{indecomposable}. 
\item Every morphism in $\C$ possesses a kernel. If $\C=\LRF$, and $$\Psi = (C_1,C_2,f):F_1 \rightarrow F_2 ,$$ then
$$ Ker(\Psi) = (C_{null}, C_1, id): P_{C_1}(F_1) \rightarrow F_1$$
If $\C=\LFG$, and 
$$ \Phi = (\gamma_1, \gamma_2, f) : \Gamma_1 \rightarrow \Gamma_2$$
then
$$ Ker(\Phi) = (\emptyset, \gamma_1, id): \gamma_1 \rightarrow \Gamma_1$$ 
\item Every morphism in $\C$ possesses a cokernel. If $\C=\LRF$, and $\Psi$ is the above morphism, then
$$ Coker(\Psi) =  (C_2, C_{full}, id): F_2 \rightarrow R_{C_2}(F_2). $$ If $\C = \LFG$, and $\Phi$ is as above, then
$$  Coker(\Phi) = (\gamma_2, \Gamma_2/ \gamma_2, id): \Gamma_2 \rightarrow \Gamma_2 / \gamma_2. $$
\item By the previous two remarks, every morphism in $\C$ has a mono-epi factorization, and moreover, the notion of exactness makes sense. 
\item \label{prop6} For a pair of objects $A, B \in \C$ such that $A \subset B$, there is an inclusion-preserving bijection between subobjects of $B$ containing $A$ and subobjects of $B/A$ such that if $\overline{S} \subset B/A$ corresponds to $S \subset B$, we have $(B/A)/\overline{S} \cong B/S$.
\item The Grothendieck group $K_{0} (\C)$ may be defined as usual (as the quotient of the free abelian group of objects by the subgroup of relations coming from short exact sequences). We have $K_0 (\LRF) = \mathbb{Z}$, while $K_{0}(\LFG)$ is the $\mathbb{Z}$--span of primitive graphs. 
\end{enumerate}

The above properties suffice to define the Ringel-Hall algebras of the categories $\LRF, \LFG$ as in section \ref{Hall_alg}. The following theorem is proven in \cite{KS}:

\begin{theorem}
  The Ringel-Hall Lie algebras of the categories $\LRF, \LFG$ are isomorphic to the  Connes-Kreimer Lie algebras on rooted trees and Feynman graphs respectively. I.e.
  $$\on{H}_{\LRF} = U(\n_{T})  \textrm{ and } \on{H}_{\LFG} = U(\n_{FG}).$$ 
\end{theorem}

In the two examples at hand, the "moduli stack" $\on{Iso}(\C)$ of isomorphism classes in $\C$ is zero-dimensional, and we may equip $\on{H}_{\C}$ with 
the symmetric bilinear form $\kappa$:
\begin{equation} \label{kappa}
\kappa(f,g) := \sum_{A \in \on{Iso}(\C)} f(A) g(A)
\end{equation}
We will use $\kappa$ to identify $\on{H}_{\C}$ with its graded dual below. The dual of $\delta_A$ with respect to $\kappa$ is denoted $\phi_A$. 

\section{Four natural representations of $\n_{\C}$} \label{left_right}

Let $\mathfrak{g}$ be a complex Lie algebra, and let $\mathfrak{g}^{op}$ denote the opposite Lie algebra, with bracket $[X,Y]_{op} := [Y,X] = - [X,Y]$. The map
\begin{align*}
{}^{t} : \mathfrak{g} & \rightarrow \mathfrak{g}^{op} \\
X & \rightarrow X^t := -X
\end{align*}
Is a Lie algebra isomorphism of $\mathfrak{g}$ with $ {\mathfrak{g}^{op}}$. 
The map ${}^t$ extends to an isomorphism 
\begin{equation} \label{transpose}
{}^t: U(\mathfrak{g}) \rightarrow U(\mathfrak{g}^{op}) \simeq U(\mathfrak{g})^{op}.
\end{equation} 
Recall that the structure of a $\mathfrak{g}$--module is equivalent to the structure of a left module over the  enveloping algebra $U(\mathfrak{g})$. The isomorphism $\ref{transpose}$
implies that a right $U(\mathfrak{g})$--module structure on $M$ is equivalent to a left $U(\mathfrak{g})$ (and hence $\mathfrak{g}$)--module structure on $M$ via the twisted action $X \cdot m := m \cdot X^t$.
To rephrase, the map ${}^t: v \rightarrow v^t$ yields an isomorphism of $\mathfrak{g}$--modules between $U(\mathfrak{\g})$ and $U(\g)^{op}$, after identifying $\mathfrak{g}$ with $\mathfrak{g}^{op}$ using ${}^t$. 
\bigskip

For the rest of the paper, $\C$ will denote one of $\LRF, \LFG$. As proven in \cite{KS}, we have $\on{H}_{\C} \simeq U(\n_{\C})$. We now consider the left and right actions of $\on{H}_{\C}$ on itself. We will write $X_{A}$ for the elment $\delta_{A}$ viewed as an element of $U(\n_{\C})$. We have
\begin{itemize}
\item {\bf Left action:}
\begin{equation} \label{ins_action}
X_A \cdot \delta_B := \delta_A \times \delta_B = \sum_{ \substack{\{ C \in \on{Iso}(\C) | \\ A \subset C, C/A \simeq B\}} }\delta_C.
\end{equation}
This is the action of $\n_{\C}$ by \emph{insertions} considered in \cite{CK}. We denote by $\on{H}_{ins}$ the vector space $\on{H}_{\C}$ with the $\n_{\C}$--module structure given by \ref{ins_action}. 
\bigskip
\begin{remark}
Let $\on{H}^{dec}_{ins}$ denote the subspace of $\on{H}_{ins}$ consisting of functions supported on decomposable objects - i.e.  $$\on{H}^{dec}_{ins} = \on{span} \{ \delta_{A} |  A \in \on{Iso}_{C}, A = B \oplus C, B \neq \emptyset, C \neq \emptyset \}. $$
This is a left ideal of $\on{H}_{\C}$ (i.e. a subrepresentation of $\on{H}_{ins}$), and the quotient $\on{H}_{ins} / \on{H}^{dec}_{ins}$ is spanned by delta functions supported on indecomposable objects. This is the insertion representation considered in \cite{CK}. 
\end{remark}

\item {\bf Right action:}
\begin{equation} \label{top_ins_action}
\delta_B \cdot X_A := \delta_B \times \delta_A = \sum_{\substack{ \{C \in \on{Iso}(\C) | \\ B \subset C, C/B \simeq A\}}} \delta_C
\end{equation}
We call this the \emph{top insertion} action. It is closely related to the operator $B_+$ considered in \cite{CK}.  We denote by $\on{H}_{top-ins}$ the vector space $\on{H}_{\C}$ with the $\n_{\C}$--module structure given by \ref{top_ins_action}. 

\begin{example} If $\C = \LRF$, we have:
\psset{levelsep=0.3cm, treesep=0.3cm}
\[
\delta_{\Tr{\bullet}} \cdot X_{ \pstree{\Tr{\bullet}}{\Tr{\bullet}\Tr{\bullet}} } = \delta_{  \pstree{\Tr{\bullet}}{\Tr{\bullet} \pstree{\Tr{\bullet}}{\Tr{\bullet}} }} +\delta_{ \Tr{\bullet} \pstree{\Tr{\bullet}}{\Tr{\bullet}\Tr{\bullet}}  } + 3 \delta_{\pstree{\Tr{\bullet}}{\Tr{\bullet}\Tr{\bullet}\Tr{\bullet} } }
\]

\end{example}

\end{itemize}

\bigskip
\noindent By the observations at the beginning of this section, $$\on{H}_{ins} \simeq \on{H}_{top-ins}$$ as $\n_{\C}$ modules. Note however, that at a combinatorial level, the isomorphism ${}^t$ is nontrivial. For instance

\begin{example}
\psset{levelsep=0.3cm, treesep=0.3cm}
\[
\delta^t_{\Tr{\bullet} \Tr{\bullet}} = \delta_{\Tr{\bullet} \Tr{\bullet}} + \delta_{\pstree{\Tr{\bullet}}{\Tr{\bullet}}}
\]
\end{example}

We may also consider the dual representations of \ref{ins_action}, \ref{top_ins_action}. Let $\on{H}^*_{\C}$ denote the restricted dual of $\on{H}_{\C}$ with basis $\{ \phi_A  | A \in \on{Iso}(\C) \}$ dual to $\{  \delta_A | A \in \on{Iso}(\C) \}$  (equivalently, $\phi_A$ is the image of $\delta_A$ under $\kappa$).

\begin{itemize}
\item The representation dual to \ref{ins_action} is determined by
\begin{align*}
X_A \cdot \phi_B (\delta_C) & = - \phi_B ( X_A \cdot \delta_C) \\ & = -\phi_B (\sum_{\substack{ \{D \in \on{Iso}(\C) | \\ A \subset D, D/A \simeq C\}}} \delta_{D}) \\  & = | \{ A \subset B | B/A \simeq C \} |
\end{align*}
which implies that
\begin{equation} \label{elim_action}
X_A \cdot \phi_B = - \sum_{A \subset B} \phi_{B/A}
\end{equation}
This is the \emph{elimination} action considered in \cite{CK} (up to a sign corresponding to an isomorphism between $\n_{\C}$ and $\n^{op}_{\C}$). We denote by $\on{H}^*_{elim}$ the vector space $\on{H}^*_{\C}$ equipped with the action \ref{elim_action}.
\item The representation dual to \ref{top_ins_action} is determined by
\begin{align*}
(\phi_B \cdot X_A) (\delta_C) &= - \phi_B (\delta_C \cdot X_A) \\ &= - \phi_B (\sum_{ \substack{\{ D \in \on{Iso}(\C) | \\ C \subset D,  D/C \simeq A \} }} \delta_D ) \\ & = - | \{ C \subset B | B/C = A \} |  
\end{align*}
which implies that
\begin{equation} \label{top_elim_action}
\phi_B \cdot X_A = - \sum_{ \substack{ \{ C \in \on{Iso}(\C) | \\ C \subset B, B/C \simeq A \} }} \phi_C
\end{equation}
We call this the \emph{top elimination} action. We denote by $\on{H}^*_{top-elim}$ the vector space $\on{H}^*_{\C}$ equipped with the action \ref{top_elim_action}.
\begin{example}
If $\C = \LRF$, we have:
\psset{levelsep=0.3cm, treesep=0.3cm}
\[
\phi_{ \pstree{\Tr{\bullet}}{\Tr{\bullet}\Tr{\bullet}} \; \; \pstree{\Tr{\bullet}}{\pstree{\Tr{\bullet}}{\Tr{\bullet}} } } \cdot X_{\pstree{\Tr{\bullet}}{\Tr{\bullet}}} = 2 \phi_{\pstree{\Tr{\bullet}} \; \; \pstree{\Tr{\bullet}}{\pstree{\Tr{\bullet}}{\Tr{\bullet}} }  } + \phi_{   \pstree{\Tr{\bullet}}{\Tr{\bullet}\Tr{\bullet}} \Tr{\bullet} }
\]

\end{example}

\end{itemize}

Again, we have $$ \on{H}^*_{elim} \simeq \on{H}^*_{top-elim}.$$

\begin{remark} \label{graded_rep}
It is a feature of Ringel-Hall algebras that $\on{H}_{\C}$ carries a natural $K_{0}(\C)$--grading. The representations constructed in this section are therefore $K_{0}(\C)$--graded representations of $\n_{\C}$. For $B \in \on{Iso}(\C)$, the gradation is given by assigning degree $[B] \in K_{0}(\C)$ to $\delta_{B}, \phi_{B}$.  
Given $A \in \on{Iso}_{\C}$, the corresponding element $X_{A} \in \n_{\C}$ has degree $[A]$ in the representations $\on{H}_{ins}, \on{H}_{top-ins}$, and degree $-[A]$ in $\on{H}^*_{elim}, \on{H}^{*}_{top-elim}$. 
\end{remark}

\section{Finite-dimensional representations}

In this section, we construct finite-dimensional sub/quotients of the four representations $\on{H}_{ins}, \on{H}_{top-ins}, \on{H}^*_{elim}, \on{H}^*_{top-elim}$ considered in the previous section.

\subsection{Finite-dimensional representations from $\on{H}^*_{elim}, \on{H}_{ins} $}
\label{fdrep1}

Let us fix an isomorphism class $M \in \on{Iso}(\C)$, and let
\[
Quot(M) := \{ E \in \on{Iso}(\C) | E \textrm{ is a quotient of } M \}
\]
(I.e. isomorphism classes of quotients of $M$). 
Let $\on{H}_{quot(M)}$ denote the subspace of $\on{H}^*_{\C}$ spanned by $\phi_A$, $A \in Quot(M)$. This is clearly a finite-dimensional subspace of $\on{H}^*_{\C}$.  It follows from property \ref{prop6} in section \ref{properties} that $\on{H}_{quot(M)}$ is stable under the action \ref{elim_action}. We thus obtain

\begin{theorem}
For every $M \in \on{Iso}(\C)$, the subspace $\on{H}_{quot(M)}$ is a finite-dimensional representation of $\n_{\C}$ with respect to the action \ref{elim_action}. 
\end{theorem}
 
If $N$ is a quotient of $M$, then any quotient of $N$ is automatically a quotient of $M$, and so in this case we have an inclusion of $\n_{\C}$--modules $\on{H}_{quot(N)} \subset \on{H}_{quot(M)}$. It is a natural question when this inclusion is split. 
Recall that an object $M$ in $\C$ is decomposable if $M \simeq A \oplus B$. In $\C = \LRF, \LFG$, $\oplus$ corresponds to disjoint union. Thus, indecomposable objects are labeled connected trees/graphs.

\begin{theorem} \label{fdtheorem}
\begin{itemize}
\item[a.)] If $M \simeq A \oplus B$, then $\on{H}_{quot(M)} = \on{H}_{quot(A)} \oplus \on{H}_{quot(B)}$ as $\n_{\C}$--modules.
\item[b.)] If $M$ is indecomposable, then $\on{H}_{quot(M)}$ is indecomposable. 
\end{itemize}
\end{theorem}

\begin{proof}
a.) Suppose $M = A \oplus B$. By the remark preceding the statement of the theorem, we have a map of $\on{H}_{\C}$--modules $\iota_A : \on{H}_{quot(A)} \hookrightarrow \on{H}_{quot(M)}$. We now construct a map $$\pi_A: \on{H}_{quot(M)} \rightarrow \on{H}_{quot(A)}.$$ Suppose that $N$ is a quotient of $M$, i.e. $N \simeq M/K$ for some $K \subset M$. Then $K = K_A \oplus K_B$, where $K_X = K \cap X, \; \; X=A,B$. Let
\[
\pi_A (\phi_N) = \phi_{A/K_A}
\]
Then $\pi_A$ is a map of $\on{H}_{\C} = U(\n_{\C})$--modules, and $\pi_A \circ \iota_A = id$. It follows that $$\on{H}_{quot(M)} = \on{H}_{quot(A)} \oplus \on{H}_{quot(B)}.$$ \\
b. ) Suppose that $ v \in \on{H}_{quot(M}$ is a nonzero vector. We claim that there exists an $X \in U(\n_{\C}) \simeq \on{H}_{\C}$ such that $ X \cdot v  = \delta_{\emptyset}$. We can write $v = \sum_{A \in \on{Iso}(\C)} c_A \phi_A \; c_A \in \mathbb{Q}$. Let $Z$ be an isomorphism class of maximal length occurring in the sum with $c_Z \neq 0$. We have $X_Z \cdot v = - c_Z \phi_{\emptyset}$. This implies that any subrepresentation of $\on{H}_{quot(M)}$ contains $\phi_{\emptyset}$, and so $\on{H}_{quot(M)}$ is indecomposable. 

\end{proof}

We can give an explicit description of the dual representation $\on{H}^*_{quot(M)}$. Using the inner product $\kappa$ to identify $\on{H}_{\C}$ with $\on{H}^*_{\C}$, and so $\on{H}^*_{quot(M)}$ with a subspace of $\on{H}_{\C}$, we see  that if $B \in Quot(M)$, then
\[
X_A \cdot \delta_B = \sum_{\{ \substack{ A \subset C | C/A = B, \\ C \in Quot(M) } \}} \delta_C
\]
I.e we insert $A$ into $B$, but sum only over those isomorphism classes which are quotients of $M$.

\subsection{Finite-dimensional representations from $\on{H}^*_{top-elim}, \on{H}_{top-ins}$}
\label{fdrep2}
Let us fix an isomorphism class $M \in \on{Iso}(\C)$, and let
\[
\leq M := \{ E \in \on{Iso}{\C} | E \textrm{ is a subobject of } M \}
\]
(I.e. isomorphism classes of subobjects of $M$). 
Let $\on{H}_{\leq M }$ denote the subspace of $\on{H}^*_{\C}$ spanned by $\phi_A,  \; A \in \leq M $. This is a finite-dimensional subspace of $\on{H}^*_{\C}$.  It follows from  the fact that subobjects of subobjects are subobjects that $\on{H}_{\leq M}$ is stable under the action \ref{top_elim_action}. We thus obtain

\begin{theorem}
For every $M \in \on{Iso}(\C)$, the subspace $\on{H}_{\leq M}$ is a finite-dimensional representation of $\n_{\C}$ with respect to the action \ref{top_elim_action}. 
\end{theorem}
 
If $N$ is a subobject of $M$, then any subobject of $N$ is automatically a subobject of $M$, and so in this case we have an inclusion of $\n_{\C}$--modules $\on{H}_{\leq N} \subset \on{H}_{\leq M}$. We have the following analogue of theorem \ref{fdtheorem}

\begin{theorem} 
\begin{itemize}
\item[a.)] If $M \simeq A \oplus B$, then $\on{H}_{\leq M} = \on{H}_{\leq A} \oplus \on{H}_{\leq B}$ as $\n_{\C}$--modules.
\item[b.)] If $M$ is indecomposable, then $\on{H}_{\leq M}$ is indecomposable. 
\end{itemize}
\end{theorem}

\begin{proof}
a.) Suppose $M = A \oplus B$. By the remark preceding the statement of the theorem, we have a map of $\on{H}_{\C}$--modules $\iota_A : \on{H}_{\leq A} \hookrightarrow \on{H}_{\leq M}$. We now construct a map $$\pi_A: \on{H}_{\leq M} \rightarrow \on{H}_{\leq A}.$$ Suppose that $N$ is a subobject of $M$. Then $N = N_A \oplus N_B$, where $N_A = N \cap A$, $N_B = N \cap B$. Let
\[
\pi_A (\phi_N) = \phi_{N_A}
\]
Then $\pi_A$ is a map of $\on{H}_{\C} = U(\n_{\C})$--modules, and $pi_A \circ \iota_A = id$. It follows that $$\on{H}_{\leq M} = \on{H}_{\leq A} \oplus \on{H}_{\leq B}.$$ \\

b. ) Suppose that $ v \in \on{H}_{\leq M}$ is a nonzero vector. We claim that there exists an $X \in U(\n_{\C}) \simeq \on{H}_{\C}$ such that $v \cdot X = \phi_{\emptyset}$. We can write $v = \sum_{A \in \on{Iso}(\C)} c_A \phi_A$. Let $Z$ be an isomorphism class of maximal length occurring in the sum with $c_Z \neq 0$. We have $v \cdot X_{Z} = - c_Z \phi_{\emptyset}$. This implies that any subrepresentation of $\on{H}_{\leq M}$ contains $\phi_{\emptyset}$, and so $\on{H}_{\leq M}$ is indecomposable. 

\end{proof}

After identifying $\on{H}^*_{\leq M}$ with a subspace of $\on{H}_{\C}$, 
the dual representation $\on{H}^*_{\leq M}$ is explicitly given by:
\[
\delta_B \cdot X_A = \sum_{\{ \substack{ B \subset C | C/B = A, \\ C \in \leq M } \}} \delta_C
\]
I.e we top-insert $A$ onto $B$, but take only those isomorphism classes which are subobjects of $M$.

\begin{remark}
 The finite-dimensional representations of $\n_{\T}$ were classified in \cite{F}, and it would be interesting to determine how the representations above fit into the classification.  
\end{remark}

\begin{remark}
Although $\on{H}^*_{elim} \simeq \on{H}^*_{top-elim}$, the subrepresentations $\on{H}_{quot(M)}$ and $\on{H}_{\leq N}$ do not correspond under the isomorphism. I.e. for a given $M \in \on{Iso}(\C)$, there is in general no $N \in \on{Iso}(C)$ such that $\on{H}_{quot(m)} \simeq \on{H}_{\leq N}$ as $\n_{\C}$--modules. 
\end{remark}

\section{Hecke correspondences and convolutions}
\label{Hecke_corr}

In this section, we show that the representations constructed above arise through the actions of Hecke correspondences in the categories $\LRF, \LFG$.  Let $\C$ denote either of these categories, and let 
\[
\heck_{\C} := \{ (A,B) | B \in \on{Iso}(\C), A \subset B \} 
\]
(i.e. the set of pairs consisting of an isomorphism class $B$ in $\C$, and a subobject $A$
of $B$). 
There are natural maps $$\pi_{i}: \heck_{\C} \rightarrow \on{Iso}(\C) \; \; i = 1,2$$
where $\pi_1 (A,B) = A$ and $\pi_2 (A,B) = B$. We also have a natural map 
\[
\pi_q:  \heck_{\C} \rightarrow \on{Iso}(\C)
\]
with $\pi_q(A,B) = B/A$. This can be summarized in the following diagram:

\begin{diagram} 
   &   &   \heck_{\C} &  &    \\
   & \ldTo<{\pi_q} & \dTo^{\pi_1} & \rdTo>{\pi_2} & \\
 \on{Iso} (\C) &   &  \on{Iso}(\C)    &   &   \on{Iso}(\C) \\  
\end{diagram}

Let $\on{H}_{\C}$ denote the Hall algebra of $\C$, i.e. the set of finitely supported $\mathbb{Q}$--valued functions on isomorphism classes of $\C$ with its associative convolution product, and $\F(\heck_{\C})$ denote the ring of $\mathbb{Q}$--valued functions on $\heck_{\C}$, with the usual commutative multiplication. 

Taking the spaces of $\mathbb{Q}$--valued functions on the sets appearing in the above diagram, we obtain the following dual diagram:

\begin{diagram} \label{dual_diag}
   &   &   \F(\heck_{\C}) &  &    \\
   & \ruTo<{\pi^*_q} & \uTo^{\pi^*_1} & \luTo>{\pi^*_2} & \\
 \on{H}_{\C} &   &  \on{H}_{\C}   &   &   \on{H}_{\C} \\  
\end{diagram}

We also have maps $\pi_{q*}, \pi_{1*}, \pi_{2*}$ in the opposite direction, defined by integrating along the fiber. I.e. 
\[
\pi_{r*} (f) (y) = \sum_{x_y \in \pi_r^{-1}(y)} f(x_y) 
\]

We can now form three different convolutions:
\begin{enumerate}
\item \begin{align*}  \on{H}_{\C} \otimes \on{H}_{\C} & \rightarrow \on{H}_{\C} \\ f \otimes g & \rightarrow \pi_{q*} (\pi^*_1(f) \pi^*_2 (g)) \end{align*}
\item 
 \begin{align*} 
\on{H}_{\C} \otimes \on{H}_{\C}  & \rightarrow  \on{H}_{\C} \\
f \otimes g & \rightarrow \pi_{2*} (\pi^*_q (f) \pi^*_1 (g))
\end{align*}
\item \begin{align*} 
\on{H}_{\C} \otimes \on{H}_{\C}  & \rightarrow  \on{H}_{\C}  \\
f \otimes g & \rightarrow \pi_{1*} (\pi^*_q (f) \pi^*_2 (g))
\end{align*}
\end{enumerate}

We proceed to analyze the action of each of them. 
Since every function on $\on{Iso}(\C)$ can be written as a linear combination of delta-functions supported on a single isomorphism class, it suffices to analyze what happens under convolution of two such delta-functions. We denote by $\delta_{A \subset B}$ the delta-function supported on the pair $A \subset B$ in $\heck_{\C}$. 

\begin{enumerate}

\item {\bf First Convolution:} For $A,B \in \on{Iso}(\C)$, we have

\begin{align} \label{calc_action1}
\pi_{q*}(\pi^*_1 (\delta_A) \pi^*_2 (\delta_B)) & = \pi_{q*} (\sum_{ \{C | A \subset C \}} \delta_{A \subset C} \times \sum_{\{ D | D \subset B \}} \delta_{D \subset B }) \\ \nonumber
&= \pi_{q*}( \sum_{A \subset B} \delta_{A \subset B}) \\ \nonumber
&= \sum_{\{ A \subset B \}} \delta_{B/A}
\end{align}

\medskip

Let us now define a right action (denoted $\bullet$) of $\on{H}_{\C}$ on itself by 
\begin{equation} \label{action1}
g \bullet f :=  \pi_{q*} (\pi^*_1(f) \pi^*_2 (g))
\end{equation}

It follows from \ref{elim_action} that after identifying $\on{H}_{\C}$ with its restricted dual $\on{H}^*_{\C}$ and switching from right to left actions, this coincides with the \emph{elimination} representation \ref{elim_action}.

\bigskip

\item {\bf Second Convolution:} For $A,B \in \on{Iso}(\C)$, we have

\begin{align} \label{calc_action2}
\pi_{2*}(\pi^*_q (\delta_A) \pi^*_1 (\delta_B)) & = \pi_{1*} (\sum_{ \{C \subset D | D/C \simeq A \}} \delta_{C \subset D} \times \sum_{\{ B \subset E \}} \delta_{B \subset E }) \\ \nonumber
&= \pi_{2*} (\sum_{\{ B \subset E | E/B = A \}} \delta_{B \subset E} ) \\ \nonumber
& = \sum_{\{ B \subset E | E/B = A \}} \delta_E
\end{align}

\medskip

We define a right action (denoted $\Box$) of $\on{H}_{\C}$ on itself by 
\begin{equation} \label{action3}
g \Box f :=  \pi_{2*} (\pi^*_q(f) \pi^*_1 (g))
\end{equation}

It follows from \ref{top_ins_action} that this coincides with the \emph{top insertion} action. 
If instead we view this as a left action of $\on{H}_{\C}$ on itself (i.e. $g$ acting on $f$ rather than vice versa), we recover the \emph{insertion} action \ref{ins_action}.

\bigskip

\item {\bf Third Convolution:} For $A,B \in \on{Iso}(\C)$, we have

\begin{align*}
\pi_{1*}(\pi^*_q (\delta_A) \pi^*_2 (\delta_B)) & = \pi_{1*} (\sum_{ \{C \subset D | D/C \simeq A \}} \delta_{C \subset D} \times \sum_{\{ E \subset B \}} \delta_{E \subset B }) \\
&= \pi_{1*} (\sum_{\{ E \subset B | B/E = A \}}) \\
& = \sum_{ \{ E \subset B | B/E = A\}} \delta_E
\end{align*}

\medskip

We define a left action (denoted $\star$) of $\on{H}_{\C}$ on itself by 
\begin{equation} \label{action2}
f \star g :=  \pi_{1*} (\pi^*_q(f) \pi^*_2 (g))
\end{equation}

It follows from \ref{top_elim_action} that after identifying $\on{H}_{\C}$ with its restricted dual $\on{H}^*_{\C}$ and switching from right to left actions, this coincides with the \emph{top elimination} representation.

\end{enumerate}

\subsection{Finite-dimensional representations}

The finite-dimensional truncations studied in sections \ref{fdrep1}, \ref{fdrep2} can also be described in the  setup of Hecke correspondences. We discuss here the construction of $\on{H}_{\leq M}$ - that of $\on{H}_{quot(M)}$ is completely analogous. 
\bigskip
\noindent For $M \in \on{Iso}(\C)$, let:
\[
\heck_{\leq M} := \{ (A,B) | B \in \on{Iso}(\C), A \subset B \subset M \}
\]
As in the previous section, we have a diagram
\begin{diagram}
   &   &   \heck_{\leq M} &  &    \\
   & \ldTo<{\pi_q} & \dTo^{\pi_1} & \rdTo>{\pi_2} & \\
 \on{Iso} (\C) &   &  \leq M    &   &   \leq M \\  
\end{diagram}
and the dual diagram
\begin{diagram}
   &   &   \F(\heck_{\leq M}) &  &    \\
   & \ruTo<{\pi^*_q} & \uTo^{\pi^*_1} & \luTo>{\pi^*_2} & \\
 \on{H}_{\C} &   &  \on{H}_{\leq M}   &   &   \on{H}_{\leq M} \\  
\end{diagram}
where $\F(\heck_{\leq M})$ denotes the space of $\mathbb{Q}$--valued functions on $\heck_{\leq M}$. 
The left action of $\on{H}_{C}$ on $\on{H}_{\leq M}$ given by
\[
f \star g := \pi_{1*} (\pi^*_q(f) \pi^*_{2} (g)) \hspace{1cm} f \in \on{H}_{\C}, g \in \on{H}_{\leq M}
\]
is the opposite of the  \emph{top elimination} representation from section \ref{fdrep2} (i.e. coincides with it after the isomorphism ${}^t:\n_{\C} \simeq \n^{op}_{\C}$), whereas the right action
\[
g \Box f :=  \pi_{2*} (\pi^*_q (f) \pi^*_{1} (g)) \hspace{1cm} f \in \on{H}_{\C}, g \in \on{H}_{\leq M}
\]
is opposite to its dual.

\end{document}